\documentclass[11pt]{amsart}
\usepackage{amscd,amssymb}
\theoremstyle{plain}
\newtheorem{Thm}[equation]{Theorem}
\newtheorem{Cor}[equation]{Corollary}
\newtheorem{Prop}[equation]{Proposition}
\newtheorem{Lem}[equation]{Lemma}
\numberwithin{equation}{section}
\title{The mass of unimodular lattices}
\author{Mikhail Belolipetsky and Wee Teck Gan}
\address{Sobolev Institute of Mathematics, Koptyuga 4, 
630090 Novosibirsk, Russia}
\address{Max Planck Institute of Mathematics, Vivatsgasse 7, 53111 Bonn, 
Germany}
\email{mbel@math.nsc.ru}

\address{Mathematics Department, Princeton University, Princeton, NJ 
08544, USA}
\address{Mathematics Department, University of California at San 
Diego, La Jolla, CA 92093, USA}
\email{wgan@math.ucsd.edu}

\begin{document}
\maketitle

\section{\bf Introduction}

An old but fundamental problem in the arithmetic theory of quadratic forms is the 
computation of the mass of a lattice $L$ in a quadratic space $(V,q)$ over a number field $F$.
Among the pioneers of its study were Smith, Minkowski and Siegel.  
After the work of Kneser, Tamagawa and Weil, this problem can be neatly formulated 
in group theoretic terms. 
More precisely, if $G$ denotes the special orthogonal group
$SO(V,q)$, then mass$(L)$ can be expressed as the volume of $G(F) \backslash G(\mathbb{A})$
with respect to a volume form $\mu_L$ associated naturally to $L$, and one can relate $\mu_L$ 
to the Tamagawa measure of $G(\mathbb{A})$ by virtue of certain local densities.  
Since the Tamagawa number of $G$ is equal to $2$, the computation of mass$(L)$ is thus reduced to
the computation of these local densities. 
\vskip 5pt

The computation of local densities is, however, not an entirely trivial task, 
especially at a $2$-adic place of $F$. 
Though there has been much work in this direction, in the introduction to his recent paper [S], 
Shimura lamented the lack of exact formulas in the literature for mass$(L)$, even for
certain special lattices. He then went on to obtain such an exact formula 
for the maximal lattice in $(V,q)$.
\vskip 5pt

The purpose of this paper is to obtain
the mass of a unimodular lattice {\em of arbitrary signature} from the 
point of view of Bruhat-Tits theory. This is achieved by
relating the local stabilizer of the lattice to a maximal
parahoric subgroup of the special orthogonal group, and appealing to 
an explicit mass formula for parahoric subgroups developed in [GHY]. 
This explicit formula is a consequence of the important
work of Prasad [P] and its extension by Gross [Gr], and can be used
to derive the results of [S] (as was done in [GHY]).
\vskip 5pt

Of course, the exact mass formula for 
positive definite unimodular lattices
is well-known (cf. for example [CS, Pg. 409]). Moreover,  
the exact formula for lattices of signature $(n,1)$ (which give rise 
to hyperbolic orbifolds) was obtained by Ratcliffe-Tschantz [RT], 
starting from the fundamental work of Siegel. 
Our approach works uniformly for unimodular lattices 
of arbitrary signature $(r,s)$ and hopefully gives a 
more conceptual way of deriving the above known results. 
The final formulas are stated in Theorems \ref{T:mass2} and \ref{T:mass3}.  

\vskip 20pt
\newpage

\section{\bf Basic Notions}

Let $F$ be a number field with ring of integers $A$.  
For each place $v$ of $F$, 
let $F_v$ denote the corresponding local field with ring of integers $A_v$ 
if $v$ is finite. 
We shall let $\mathbb{A}$ be the adele ring of $F$ and 
$\mathbb{A}_f$ the ring of finite adeles.

\vskip 10pt

Let $(V,q)$ be a quadratic space over $F$ of dimension $d$, and
let $B_q$ be the symmetric bilinear form defined by:
\[ B_q(x,y) = \frac{1}{2} (q(x+y) - q(x) -q(y)), \quad x,y \in V. \]
Note that
\[ B_q(x,x) = q(x). \]
The discriminant of $(V,q)$ is
defined as follows. One diagonalizes $q$ using a suitable basis of $V$, 
say $q(x_1,...,x_d) = \sum_i a_i x_i^2$; then one sets
\[ disc(q) = (-1)^{d(d-1)/2} \cdot  \prod_i a_i \in F^{\times}/F^{\times 2}.
\]
The square class of $disc(q)$ determines a quadratic character
$\chi_{disc(q)}$ of $Gal(\overline{F}/F)$, and we let
$E_{disc(q)}$ be the \'{e}tale quadratic algebra determined by $\chi_{disc(q)}$. 
More precisely,
\[ E_{disc(q)} = \begin{cases}
F \times F, \text{   if $disc(q) \in F^{\times 2}$;} \\
F(\sqrt{disc(q)}), \text{   if $disc(q) \notin F^{\times 2}$.} \end{cases} \]
Further, if $E$ is an \'{e}tale quadratic algebra over $F$ as above, we let
\[ d_E = \# A/disc_{E/F} = |\mathbb{N}_{F/\mathbb{Q}}(disc_{E/F})| \]
where $disc_{E/F}$ is the discriminant ideal of the ring of integers $A_E$ over $A$.
For example, if $E = F \times F$, then $d_E  =1$. 
\vskip 5pt

Let $L \subset V$ be a lattice on which $q$ takes integer values. 
We shall write $L_v$ for the localization $L \otimes A_v$.
Recall the following basic definition:
\vskip 5pt

\noindent{\bf Definition:} 
The {\em genus} of $L$ is the set of 
isomorphism classes of lattices $M$ in $(V,q)$ such that
\[ (L_v, q) \cong (M_v, q) \quad \text{  for all finite places $v$.}  \]
\vskip 5pt

Let $G = SO(V,q)$ be the special orthogonal group associated with $(V,q)$.
To avoid having to work with non-semisimple groups, we assume henceforth that 
$d =dim(V) \geq 3$. This is for uniformity of exposition and is not a serious
assumption. If 
\[ K_L = \prod_{v < \infty} K_{L_v} \subset G(\mathbb{A}_f) \]
is the stabilizer of $\widehat{L} = \prod_{v < \infty} L_v$,
then the genus of $L$ can be indexed by the double coset space
\[ G(F) \backslash G(\mathbb{A}_f)/ K_L. \]
This is a finite set and for $\alpha \in
G(F) \backslash G(\mathbb{A}_f)/K_L$, the corresponding lattice in $(V,q)$
is given by
\[ L_{\alpha} = V \cap \alpha \widehat{L}, \]
with the intersection occurring in $V \otimes \mathbb{A}_f$. 
For each $\alpha$, let
\[ \Gamma_{\alpha} = G(F) \cap \alpha K_L \alpha^{-1}. \]
It is the stabilizer in $G(F)$ of $L_{\alpha}$ and 
has finite covolume in $G(F \otimes \mathbb{R})$.
\vskip 10pt

We can now introduce an important classical invariant for the genus of $L$;
it is called the mass of $L$.
When $F$ is totally real and $q$ is totally definite, this is defined by
\[ Mass(L)= \sum_{\alpha} \frac{1}{\# \Gamma_{\alpha}}. \]
Another way to define this is as follows. Let $\mu_L$ be the Haar measure on 
$G(\mathbb{A})$ which gives the open compact subgroup
$G(F \otimes \mathbb{R}) \times K_L$ volume 1. Then
\[ Mass(L) = \int_{G(F) \backslash G(\mathbb{A})} \mu_L. \]
\vskip  5pt

Using the above integral formula, one can extend the definition of $Mass(L)$
to the indefinite case. However, to 
define the measure $\mu_L$ in general, one needs to 
specify a Haar measure on $G(F \otimes \mathbb{R})$. Different authors made
different choices for this; we shall now explain our choice.
\vskip 5pt

There are two natural choices of Haar measure on $G(F \otimes \mathbb{R})$, 
at least from the group theoretic 
point of view. Any reductive algebraic group 
over $\mathbb{R}$, such as $Res_{F/\mathbb{Q}}\;G \times \mathbb{R}$, 
has a unique split form and a unique compact form
and each of these has a natural Haar measure. 
For the split form, it is the
measure induced by an invariant differential form of 
top degree on the canonical Chevalley model over $\mathbb{Z}$. 
For the compact form, it is the measure
giving the group volume 1. Each of these two measures can be transferred 
to any other forms of the group in a standard way, as described in [Gr, \S 11].
In this way, we obtain two natural Haar measures 
$\mu_{cpt}$ and $\mu_{sp}$ on $G(F \otimes \mathbb{R})$. They are related by (cf. [Gr, \S 7])
\[ \mu_{sp} = \gamma_G^{deg(F)} \cdot \mu_{cpt} \]
where $deg(F)$ is the degree of $F$ over $\mathbb{Q}$ and
\[ \gamma_G = \begin{cases}
\frac{(2\pi)^{n(n+1)}}{\prod_{r=1}^{n}(2r-1)!}, \text{  if $d = 2n+1$;} \\
\frac{(2\pi)^{n^2}}{(n-1)! \cdot 
\prod_{r=1}^{n-1} (2r-1)!}, \text{  if $d = 2n$.} 
\end{cases} \]

\vskip 5pt

The measure we use for the definition of $Mass(L)$ is $\mu_{cpt}$.
Thus we set
\[ \mu_L = \mu_{cpt} \times \prod_{v < \infty} \mu_{L_v} \]
where $\mu_{L_v}$ is the Haar measure of $G(F_v)$ giving $K_{L_v}$ volume 1.
Then $Mass(L)$ is defined in general by:
\[ Mass(L) = \int_{G(F) \backslash G(\mathbb{A})} \mu_L. \]
When $q$ is totally definite, 
this agrees with  the classical definition above.
When $q$ is indefinite, we have:
\[ Mass(L) = \sum_{\alpha} \int_{\Gamma_{\alpha} \backslash 
G(F \otimes \mathbb{R})} \mu_{cpt}. \]

\vskip 10pt

\section{\bf A Mass Formula}
For a natural class of lattices $L$ in $(V,q)$, an explicit formula for 
$Mass(L)$ was given in [GHY, Proposition 2.13], based on the fundamental 
work [P] and its extension [Gr]. These are the lattices for 
which the local stabilizers $K_{L_v}$ are parahoric subgroups of 
$G(F_v)$. To state this formula, we need to introduce some more notations.

\vskip 5pt

As in [GHY], our usage of the word ``parahoric'' is slightly different
from that in Bruhat-Tits theory. Namely, we call an open compact 
subgroup $K_{L_v}$ ``parahoric'' if it stabilizes setwise a simplex in the 
Bruhat-Tits building  of $G(F_v)$.
By Bruhat-Tits theory, there is a smooth affine group scheme 
$\underline{G}_v$ over $A_v$, with generic fiber $G \times_F F_v$, such that
\[ K_{L_v} = \underline{G}_v(A_v). \]
With our notion of parahoric subgroups, the group scheme $\underline{G}_v$ 
is possibly disconnected here. Let $\overline{G}_v$ be the 
maximal reductive quotient of the special fiber of $\underline{G}_v$;
it is a (possibly disconnected) reductive algebraic group over the residue 
field of $F_v$. Let $N(\overline{G}_v)$ be the number of positive 
roots of $\overline{G}_v$ over the algebraic closure. Then the number 
\[ q_v^{-N(\overline{G}_v)} \cdot \# \overline{G}_v(A_v/\pi_v) \]
is a product of certain local $L$-factors and can be easily computed once one 
identifies $\overline{G}_v$.
\vskip 5pt

Similarly, if $G_{qs}$ is the quasi-split inner form of $G$, then
we may consider the connected integral model 
$\underline{G}_{qs,v}$ associated 
to the special maximal compact subgroup of $G_{qs}(F_v)$ specified
in [Gr, \S 4].  As above, we may define the number
\[ q_v^{-N(\overline{G}_{qs,v})} \cdot \# \overline{G}_{qs,v}(A_v/\pi_v). \]
We can now state the formula of [GHY, Proposition 2.13]:

\begin{Thm} \label{T:mass}
Let $L$ be a lattice in $(V,q)$ so that the local stabilizer $K_{L_v}$ 
is a parahoric subgroup of $G(F_v)$ for each finite place $v$. 
Then
\[ Mass(L) = \left( \tau(G) \cdot d_F^{dim(G)/2} \cdot
\gamma_G^{-deg(F)}\cdot L(G) 
\right) \cdot \prod_v \lambda_{L_v}, \]
where
\begin{itemize}
\item[-] $\tau(G) =2$ is the Tamagawa number of $G$;

\item[-]  $d_F$ is the absolute value of the discriminant of $F/\mathbb{Q}$;

\item[-] $L(G)$ is the special value of an $L$-function associated to $G$ and is
given by:
\[ L(G) = \begin{cases}
\prod_{r=1}^n \zeta_F(2r) \text{   if $d = 2n+1$;} \\
\prod_{r=1}^{n-1} \zeta_F(2r) \cdot L(n, \chi_{disc(q)}) \cdot d_{E_{disc(q)}}^{n-\frac{1}{2}}, \text{  if
$d = 2n$.} \end{cases} \]

\item[-] for each finite place $v$,
\[ \lambda_{L_v} = 
\frac{q_v^{-N(\overline{G}_{qs,v})} 
\cdot \# \overline{G}_{qs,v}(A_v/\pi_v)}
{q_v^{-N(\overline{G}_v)} \cdot \# 
\overline{G}_v(A_v/\pi_v)}. \]
\end{itemize}
\end{Thm}
\vskip 10pt

\noindent{\bf Remarks:}
(i) The first factor in the mass formula depends only on the group $G$ and 
should be regarded as the main term,  whereas the $\lambda$-factors
depend on the local stabilizers $K_{L_v}$ 
and should be regarded as fudge factors.
The point of the above formula is that the 
$\lambda$-factors are effectively computable from Bruhat-Tits theory. 
Note that the product of $\lambda$-factors is a finite product, since 
for almost all places $v$, $G_v$ is quasi-split and  
$K_{L_v}$ is a hyperspecial maximal compact subgroup, in which case 
$\lambda_{L_v}= 1$. 
\vskip 5pt

(ii) In [GHY, Proposition 2.13], the formula was stated only under the assumption
that $F$ is totally real and $q$ is totally definite. However, the derivation
given in [GHY, \S2] only relies on [GG, Proposition 9.3], and the latter  
holds without these restrictions, as long as the group $G$ is semisimple.
\vskip 5pt

(iii) The reader may notice that the formula in [GHY, Proposition 2.13] 
is cleaner than the one given above, and involves $L$-values at 
negative integers rather than positive integers. Not surprisingly,
the two versions are related by the functional equation. 
If we had applied the functional equation to the formula of the Theorem,
we expect to get the values $\zeta_F(1-2r)$. But for some number fields, 
these quantities are equal to $0$; in this case, one has to use the leading 
term of the Taylor expansion of the zeta functions in place of the value.
Because of this complication, we prefer to leave the formula as it is.
\vskip 15pt

Suppose that $L$ is an arbitrary lattice, so that $K_{L_v}$ may not be
parahoric. For each $v$, there will be a maximal 
parahoric subgroup $K_v$ containing
$K_{L_v}$, and clearly, if we know the index of $K_{L_v}$ in $K_v$, 
we can determine $Mass(L)$. More precisely, if
we set 
\[ \lambda_{L_v} = \#K_v/K_{L_v} \cdot \lambda_{K_v},\] 
then $\lambda_{L_v}$ is well-defined, i.e. independent 
of the choice of $K_v$,
and the formula of the Theorem continues to hold for any lattice $L$.
\vskip 10pt

In the following sections, we shall use this formula to obtain the mass
of certain special lattices $L$.

\vskip 10pt

\section{\bf Maximal and Unimodular Lattices}

We recall the following basic definitions:

\vskip 5pt

\noindent{\bf Definitions:} 
\begin{itemize}
\item  $L$ is a {\em maximal lattice} in $(V,q)$ if $q$  
is not $A$-valued on any lattice strictly containing $L$.

\item $L$ is a {\em unimodular lattice}
 in $(V,q)$ if $L$ is self-dual with respect to $B_q$,
i.e. $L^* = L$, where
\[ L^* = \{ x \in V:  B_q(x,L) \subset A \}. \]
\end{itemize}
Clearly, the analogous 
definitions can be made for the local lattices $L_v$. We note the following elementary remarks.
\vskip 5pt

\noindent{\bf Remarks:} 
(i) $L$ is maximal if and only if 
$L_v$ is maximal for all finite places $v$ of $F$.
Similarly, $L$ is unimodular if and only if $L_v$ 
is unimodular for all finite places $v$.

\vskip 5pt

(ii) If $L$ is a unimodular lattice, then of course $L$ is 
a largest lattice on which $B_q$ is $A$-valued. If $F_v$ is a $p$-adic field 
with $p \ne 2$, then this implies that $L_v$ is a maximal lattice.
However, over the number field $F$ or a 2-adic field $F_v$, $L$ need not be 
a maximal lattice. 

\vskip 5pt

For a maximal lattice $L$, it is a consequence of [BT] that 
the local stabilizers $K_{L_v}$ are maximal parahoric subgroups.
Exploiting this fact, the
mass of $L$ can be obtained using Theorem \ref{T:mass}. This was carried out
in [GHY], where the relevant $\lambda$-factors were tabulated.
\vskip 5pt

In this paper, we shall explain how to obtain the mass of a 
unimodular lattice by exploiting Theorem \ref{T:mass}, at least for a quadratic 
space over $\mathbb{Q}$ (in fact, over any number field $F$ such that 
for any place $v$ lying over the prime $2$, $F_v$ is unramified over 
$\mathbb{Q}_2$). To do this, we need to relate the local
stabilizer of a unimodular lattice to a parahoric subgroup.  
We treat this local question in the next section. 
\vskip 10pt

\section{\bf Local stabilizers of Unimodular Lattices}
In this section, 
we shall relate the stabilizers of unimodular lattices to
maximal parahoric subgroups. Recall that we are assuming that $d =dim(V) \geq 3$.
For simplicity, we assume that $F = \mathbb{Q}_p$. However,
our discussion holds for any finite extension of 
$\mathbb{Q}_p$ if $p$ is odd, and can be extended to cover any unramified 
finite extension of $\mathbb{Q}_2$. 

\vskip 10pt

First recall the classification of quadratic spaces over $\mathbb{Q}_p$.
We have already defined the discriminant of $(V,q)$.
Another invariant of $(V,q)$ is 
the Hasse-Witt invariant defined as follows.
By choosing a suitable
basis of $V$, we may diagonalize the form $q$, say
\[ q(\underline{x}) = \sum_{i=1}^d a_i x_i^2. \] 
Then the Hasse-Witt invariant is:
\[ 
\epsilon_{HW}(V,q) = \prod_{i < j} (a_i,a_j) \in \{ \pm 1 \}, \]
where $(-,-)$ is the Hilbert symbol of $\mathbb{Q}_p$.
The quadratic space $(V,q)$ is then determined by the invariants
\[ (dim(V), disc(V,q), \epsilon_{HW}(V,q)). \]
Note that the definition of $\epsilon_{HW}$ differs from the $\epsilon$ in 
[GHY]. 

\vskip 10pt

As we mentioned in the previous section, if $p\ne 2$, a unimodular lattice 
in $(V,q)$ is necessarily a maximal lattice. These were enumerated in [GHY],
with their stabilizers identified and their $\lambda$-factors tabulated.
Hence the main local problem is to understand the stabilizers of unimodular lattices
over $\mathbb{Q}_2$. 
\[ \text{{\em We shall assume for the rest of the section that
$p =2$}.} \] 
In this case, one distinguishes betwen two types of unimodular 
lattices.
\vskip 10pt

\noindent{\bf Definitions:}
A unimodular lattice $L$ in $(V,q)$ is said to be {\em even} if 
$q(L) \subset 2\mathbb{Z}_2$; it is said to be {\em odd} otherwise.
\vskip 10pt

\begin{Prop} \label{P:uni}
Fix a quadratic space $(V,q)$. Then 
there is at most one isomorphism class of 
odd (resp. even) unimodular lattices in $(V,q)$. 
\end{Prop}

\vskip 5pt

\begin{proof}
This follows from [OM, Theorem 93.16, Pg. 259] and the discussion in 
[OM, \S 93G].
\end{proof}

To decide which quadratic spaces actually possess a unimodular lattice, 
we first begin with the even case. We have:

\begin{Prop} \label{P:even1}
(i) If $L$ is an even unimodular lattice in $(V,q)$, then
$L$ is a maximal lattice  in $(V, \frac{1}{2}q)$. 
Moreover, $d=2n$ is even and 
$disc(q)$ can be represented by an element of $\mathbb{Z}_2^{\times}$.
\vskip 5pt

(ii) $(V,q)$ contains an even unimodular lattice if and only if $E_{disc(q)}$
is not a ramified quadratic extension and $\epsilon_{HW}(q) = 
(-1)^{n(n-1)/2}$. 
\vskip 5pt

(iii) More explicitly, the possible quadratic spaces are:
\begin{itemize}
\item $(V,q) \cong \mathbb{H}^n$, where $\mathbb{H} = \langle e, f \rangle$
is the hyperbolic plane. An even unimodular lattice is
\[ L_{even} = \langle e_1,...,e_n, 2f_1,...,2f_n \rangle. \]

\item $(V,q) = (E, 2 \mathbb{N}_E) \oplus \mathbb{H}^{n-1}$, where $E$ is the
unramified quadratic extension of $\mathbb{Q}_2$ with norm map $\mathbb{N}_E$.
An even unimodular lattice is
\[ L_{even} = A_E \oplus \langle e_i,...,e_{n-1}, 2f_1,...,2f_{n-1} \rangle. 
\] 
\end{itemize}
\end{Prop}

\begin{proof}
The first assertion of (i) is clear.
The fact that $d$ is even was shown in [OM, 93.15, Pg. 258] and 
the statement about $disc(q)$ is obvious.
\vskip 5pt

By (i), to obtain an explicit list of $(V,q)$ which contains
even unimodular lattices, it suffices to
examine a maximal lattice in $(V,\frac{1}{2}q)$ and see if 
it is self-dual with respect to $B_q$. Using the enumeration
of maximal lattices in [GHY], a short check gives the list in (iii)
and it is easy to see that these two quadratic spaces have 
the discriminant and Hasse-Witt invariant stated in (ii).
\end{proof}

\begin{Cor}
If $L$ is even unimodular, then 
$K_L$ is the stabilizer in $G(\mathbb{Q}_2) 
= SO(V,q)(\mathbb{Q}_2) = SO(V,\frac{1}{2}q)(\mathbb{Q}_2)$
of a maximal lattice in $(V ,\frac{1}{2}q)$. In particular, it
is a maximal parahoric subgroup and  $\lambda_L =1$.
\end{Cor}

Now we come to the odd unimodular lattices; the situation 
here is more interesting. It is  
not difficult to enumerate the odd unimodular lattices in the spirit 
of the previous proposition. However, we shall refrain from doing so at 
the moment, since we would like to avoid case-by-case analysis as much as 
possible. We begin by noting:

\begin{Lem} \label{L:odd1}
$(V,q)$ contains an odd unimodular lattice if and only if
$disc(V,q)$ can be represented by an element of $\mathbb{Z}_2^{\times}$.
\end{Lem}

\begin{proof} 
The ``only if'' part is clear. Conversely,
for given $\delta \in \mathbb{Z}_2^{\times}$
and $\epsilon = \pm 1$, consider the lattice $L_{\delta,\epsilon}$ defined
by the following quadratic form on $\mathbb{Z}_2^d$: 
\[ 
q_{\delta,\epsilon}(x) 
= (-1)^{[d/2]} \delta x_1^2 + \epsilon x_2^2 + 
\epsilon x_3^2 +\sum_{i= 4}^d x_i^2. \]
This defines a unimodular lattice and 
a quick check shows that 
\[ disc(q_{\delta, \epsilon}) = \delta \quad \text{and}
\quad \epsilon_{HW}(q_{\delta, \epsilon}) = \epsilon. \]
This proves the reverse implication.
\end{proof}

Let $L$ be an odd unimodular lattice in $(V,q)$.
The rest of the section is devoted to the determination of $K_L$ and
the computation of $\lambda_L$. A simple but crucial 
observation is that the induced map
\[ \bar{q}: L \rightarrow \mathbb{Z}/2\mathbb{Z} \]
is a group homomorphism, which is surjective since $L$ is odd. 
Let $\Lambda$ be the kernel of $\bar{q}$; 
it is a sublattice with index 2 in $L$. 
\vskip 5pt

We want to relate $K_L$ to $K_{\Lambda}$.
From the definition of $\Lambda$, the following lemma is clear.

\begin{Lem} \label{L:stab}
We have: $K_L \subset K_{\Lambda}$.
\end{Lem}

\vskip 5pt
Now the quadratic form $\frac{1}{2}q$  takes integer value on $\Lambda$. 
So we may ask if $\Lambda$ is a maximal lattice in $(V,\frac{1}{2}q)$.
We have:

\begin{Prop} \label{P:odd2}
Let $L$ be an odd unimodular lattice in $(V,q)$ and let $\Lambda$ be 
defined as above.
\vskip 5pt

(i) If $(V,q)$ does not contain an even unimodular lattice,
then $\Lambda$ is a maximal lattice in $(V,\frac{1}{2}q)$. 
\vskip 5pt

(ii) If $(V,q)$ contains an even unimodular lattice, then
$\Lambda$ is not a maximal lattice in $(V, \frac{1}{2}q)$.
There is an even unimodular lattice 
$L_{even}$ of $(V,q)$ such that $\Lambda = L \cap L_{even}$.
Moreover, $[L_{even}: \Lambda] =2$.
\end{Prop}

\begin{proof}
(i) Suppose that $\Lambda' \supset \Lambda$ and 
$\frac{1}{2}q$ is integer-valued 
on $\Lambda'$. Then the symmetric bilinear form $B_q$ is integer valued on 
$\Lambda'$. 
Since $\Lambda$ is contained in a self-dual lattice with index $2$,
this forces $\Lambda'$ to be self-dual with respect to $B_q$ as well.
So $\Lambda'$ is an even unimodular lattice in $(V,q)$. 
But $(V,q)$ does not contain such a lattice by assumption, and so
(i) is proved. 

\vskip 5pt

(ii) The two quadratic spaces 
listed in Proposition \ref{P:even1}(iii) do contain odd unimodular lattices.  
Indeed, an odd unimodular lattice
in $\mathbb{H}^n$ is:
\[ L_{odd} = \langle e_1,...,e_{n-1},
2f_1,...,2f_{n-1} \rangle \oplus \langle e_n+f_n, e_n-f_n \rangle \]
and one for $(E, 2\mathbb{N}_E) \oplus \mathbb{H}^{n-1}$ is
\[ L_{odd} = A_E \oplus \langle e_2,...,e_{n-1},
2f_2,...,2f_{n-1} \rangle \oplus \langle e_n+f_n, e_n-f_n \rangle. \]
If $L_{even}$ is the lattice defined in Proposition \ref{P:even1}(iii), then
one sees easily that $\Lambda = L_{even} \cap L_{odd}$.
This proves (ii).
\end{proof}
\vskip 10pt

To determine the index of $K_L$ in $K_{\Lambda}$, we observe that
\[ \Lambda \subset_2 L \subset_2 \Lambda^* \]
where $\Lambda^*$ is dual of $\Lambda$ with respect to $B_q$.
The following lemma determines the order 4 group $\Lambda^*/\Lambda$:

\begin{Lem} \label{L:quot}
\[ \Lambda^*/\Lambda = \begin{cases}
\mathbb{Z}/4\mathbb{Z}, \text{  if $d$ is odd;} \\
\mathbb{Z}/2\mathbb{Z} \times \mathbb{Z}/2\mathbb{Z}, 
\text{  if $d$ is even.} \end{cases} \]
\end{Lem}

\begin{proof}
Let $L_{\delta,\epsilon}$ be as given in the proof 
of Lemma \ref{L:odd1}, so that $L_{\delta,\epsilon} = \mathbb{Z}_2^d$ and
$q$ has the form:
\[ q(x) = \sum_i a_i x_i^2, \]
where each $a_i$ is a unit in $\mathbb{Z}_2$.
Since $a_i \equiv 1$ (mod 2), we see that
\[ \Lambda = \{ x \in L: \sum_i x_i \equiv 0 \text{(mod 2)} \} \]
and 
\[ \Lambda^* = \langle L, x_0 = \frac{1}{2} \sum_i e_i \rangle. \]
Now $\Lambda^*/\Lambda \cong \mathbb{Z}/4\mathbb{Z}$ if and only if
$x_0$ has order 4 in $\Lambda^*/\Lambda$, and this occurs if and only if
$2x_0$ does not lie in $\Lambda$. 
Clearly, this holds if and only if $d$ is odd.
\end{proof}
\vskip 5pt

Now when $d$ is odd, $L/\Lambda$ can be characterized as the unique 
order 2 subgroup of $\Lambda^*/\Lambda \cong \mathbb{Z}/4\mathbb{Z}$.
This implies that $K_{\Lambda} \subset K_L$, and together with  
Lemma \ref{L:stab}, we deduce that $K_L  = K_{\Lambda}$, which is a maximal 
parohoric subgroup of $G(F_v)$. Hence, $\lambda_L = \lambda_{\Lambda}$
can be read off from the appropriate table in [GHY]. We have shown:

\begin{Thm} \label{T:odd3}
Let $L$ be an odd unimodular lattice in $(V,q)$.
If $d=2n+1$ is odd, then $K_L$ is the stabilizer of a maximal lattice in
$(V, \frac{1}{2}q)$. Further, the value of $\lambda_L$ is given by the 
following table.
\vskip 10pt

\begin{center}
\begin{tabular}{|c|c|c|}
\hline
& & \\
$(V,q)$ & split & non-split \\
& & \\
\hline
& & \\
$\lambda_L$ & $(2^n+1)/2$ & $(2^n-1)/2$ \\
& & \\
\hline
\end{tabular}
\end{center}
\end{Thm}
\vskip 10pt

Assume henceforth that $d=2n$ is even.
In this case, Lemma \ref{L:quot} is not good enough to pinpoint
$K_L$; we need more structures on $\Lambda^*/\Lambda$.  
Indeed, the quadratic form $q$ induces a quadratic form on the
$\mathbb{F}_2$-vector space $\Lambda^*/\Lambda$ and we want to identify
this quadratic space.  The non-trivial elements 
of $\Lambda^*/\Lambda$ can be represented  by 
\[ 
x_0  = \frac{1}{2}\sum_i e_i, \quad  y_0 = e_1 \quad \text{and} \quad
z_0 = x_0-y_0. \]
Further,
\[ q(x_0) = \frac{1}{4}\sum_i a_i = q(z_0) \quad \text{and} \quad
q(y_0) =a_1. \]
Working with $q= q_{\delta,\epsilon}$, we have
\[ \sum_i a_i = (-1)^{[d/2]}\delta + 2 \epsilon(q) + d-3   \]
where $\delta \in \mathbb{Z}_2^{\times}$ is a representative of $disc(q)$.  
Clearly, the isomorphism class of the quadratic space $\Lambda^*/\Lambda$
depends on the valuation of the element
\[ \kappa(q) = (-1)^{[d/2]} \delta + 2\epsilon(q) + d-3. \]
We have the following 3 cases:
\vskip 5pt

\begin{itemize}
\item[{\bf Case A:}] ord$(\kappa(q)) =1$. In this case, 
$\frac{1}{2}q$ induces a quadratic form on $\Lambda^*/\Lambda$. 
With respect to the basis $\{ x_0, z_0 \}$, it is given by $x^2+z^2$.

\vskip 5pt

\item[{\bf Case B:}] ord$(\kappa(q)) =2$. In this case, $q$ induces a 
quadratic form on $\Lambda^*/\Lambda$ isomorphic to 
$(\mathbb{F}_4, \mathbb{N}_{\mathbb{F}_4})$.
\vskip 5pt

\item[{\bf Case C:}] ord$(\kappa(q)) \geq 3$. 
In this case, $q$ induces a quadratic form on $\Lambda^*/\Lambda$. 
With respect to the basis $\{ x_0, z_0 \}$, it is given by $xz$.
In other words, $\Lambda^*/\Lambda$ is a hyperbolic plane.
\end{itemize}

\vskip 5pt

\begin{Prop} \label{P:odd4}
(i) Cases A occurs if and only if $E_{disc(q)}$ is a ramified quadratic 
extension of $\mathbb{Q}_2$. In this case, $K_L = K_{\Lambda}$ 
and so $\lambda_L = \lambda_{\Lambda}$.
\vskip 5pt

(ii) Case B occurs if and only if $E_{disc(q)}$ is not ramified and
$\epsilon_{HW}(q) = -(-1)^{n(n-1)/2}$. In this case, 
$K_L$ has index 3 in $K_{\Lambda}$ and so $\lambda_L =
3\lambda_{\Lambda}$.
\vskip 5pt

(iii) Case C occurs if and only if $E_{disc(q)}$ is not ramified and 
$\epsilon_{HW}(q) = (-1)^{n(n-1)/2}$, i.e. $(V,q)$ contains even 
unimodular lattices. In this case, $K_L = K_{\Lambda}$ and so 
$\lambda_L = \lambda_{\Lambda}$.
\end{Prop}
\vskip 5pt

\begin{proof}
The characterization of the various cases in terms of discriminant and 
Hasse-Witt invariant is a straightforward check; we omit the details.
\vskip 5pt

In Case A, $L/\Lambda$ is the unique isotropic line in the quadratic space 
$\Lambda^*/\Lambda$; so any element of  $K_{\Lambda}$ has to fix 
$L/\Lambda$. In Case C, $L/\Lambda$ is the unique non-isotropic line
in $\Lambda^*/\Lambda$, and so again $K_{\Lambda} \subset K_L$. 
\vskip 5pt

In Case B, each of the 3 lines in $\Lambda^*/\Lambda$ is non-isotropic
and gives rise to an odd unimodular lattice. We need to show that $K_{\Lambda}$
acts transitively on these lines. If $L_1/\Lambda$ and $L_2/\Lambda$
are 2 such lines, then by Proposition \ref{P:uni}, there is an element 
$g \in G(\mathbb{Q}_2)$ such that $g(L_1) = L_2$. But $\Lambda$ can be 
characterized as the subset of $L_1$ (resp. $L_2$) on which $q$ takes 
even-integer values. So we must have $g(\Lambda) = \Lambda$. In other words, 
we have found an element of $K_{\Lambda}$ which takes $L_1$ to $L_2$.
\end{proof}
\vskip 10pt

The proposition allows us to compute $\lambda_L$ in Cases A and B,
since $\Lambda$ is a maximal lattice in $(V, \frac{1}{2}q)$ in these  cases.
The values of $\lambda_L$ are tabulated at the end of this section.
In Case C, knowing that $\lambda_L = \lambda_{\Lambda}$ does
not help us since $\Lambda$ is not a maximal lattice; 
we need to do some more work.
\vskip 10pt

Assume hence that $(V,q)$ contains even unimodular lattices. 
As we saw in Proposition \ref{P:odd2}(ii), there is an even unimodular lattice 
$L_{even}$ such that
\[ \Lambda= L_{even} \cap L. \]
Now $\frac{1}{2}q$ induces a quadratic form on $L_{even}/2L_{even}$, and
using the two $L_{even}$'s given in Proposition \ref{P:even1}(iii), it is easy to 
check that 
\[ L_{even}/2L_{even} \cong \begin{cases}
\mathbb{H}^n, \text{   if $(V,q)$ is split;} \\
\mathbb{F}_4 \oplus \mathbb{H}^{n-1}, \text{  otherwise.} \end{cases} \]
Hence,
\[ SO(L_{even}/2L_{even}) \cong  \begin{cases}
SO_{2n}, \text{  if $(V,q)$ is split;} \\
^2SO_{2n}, \text{  otherwise.} \end{cases} \]
Moreover, the action of $K_{L_{even}}$ on $L_{even}/2L_{even}$ 
gives a surjection 
\[ r : K_{L_{even}} \longrightarrow SO(L_{even}/2L_{even}). \]
Now we note:

\begin{Prop} \label{P:odd5}
(i) $K_{\Lambda} \cap K_{L_{even}}$ has index 2 in $K_{\Lambda}$.
\vskip 5pt

(ii) $K_{\Lambda} \cap K_{L_{even}}$  is the subgroup of $K_{L_{even}}$
stabilizing an isotropic line in $L_{even}/2L_{even}$.
Its index in $K_{L_{even}}$ is given by:
\[ \begin{cases}
(2^{n-1}+1)(2^n-1), \text{  if $(V,q)$ is split;} \\
(2^{n-1}-1)(2^n+1), \text{  otherwise.} \end{cases} \]
\end{Prop}
\vskip 5pt

\begin{proof}
(i) This is because $L_{even}/\Lambda$ is one of the 
two isotropic lines in $\Lambda^*/\Lambda$, 
and $K_{\Lambda}$ acts trasnsitively on the set of isotropic lines.
\vskip 5pt

(ii) The first statement follows since
$2\Lambda^*/2L_{even}$ is an isotropic line in the quadratic space
$L_{even}/2L_{even}$. Further, $SO(L_{even}/2L_{even})$ 
acts transitively on the set of isotropic lines, and the stabilizer of
one such line is a maximal parabolic subgroup $P$ with Levi factor 
\[ \begin{cases}
GL_1 \times SO_{2n-2}, \text{   if $(V,q)$ is split;} \\
GL_1 \times {^2SO_{2n-2}}, \text{   otherwise.} \end{cases} \]
Hence $K_{\Lambda} \cap K_{L_{even}}$ is equal to the 
inverse image of $P(\mathbb{F}_2)$ under the projection $r$
(which is a non-maximal parahoric subgroup). Thus
\[  \# K_{L_{even}}/K_{\Lambda} \cap K_{L_{even}}
= \#SO(L_{even}/2L_{even})/P(\mathbb{F}_2), \]
and one obtains the values listed in (ii).
\end{proof} 
\vskip 10pt

The proposition allows us to compute volume of $K_L = K_{\Lambda}$
given the volume of $K_{L_{even}}$. Thus we can compute the value of 
$\lambda_L = \lambda_{\Lambda}$ from the value of $\lambda_{L_{even}}$.
The latter is known from [GHY] since
$L_{even}$ is a maximal lattice in $(V, \frac{1}{2}q)$; 
in fact, $\lambda_{L_{even}} =1$. 
\vskip 10pt

\noindent{\bf Remarks:} Using the results of [BT], one can check that
$K_{\Lambda}$ is a maximal compact subgroup. It has an associated integral 
group scheme $\underline{G}_{\Lambda}$ whose special fiber 
has maximal reductive quotient 
\[ \overline{G}_{\Lambda} = \begin{cases}
\text{$S(O_2 \times O_{2n-2})$ if $(V,q)$ is split;} \\
\text{$S(O_2 \times$ $^2O_{2n-2})$ otherwise.} 
\end{cases} \]
\vskip 10pt

To conclude this section, we tabulate the values of $\lambda_L$ when 
$d=2n$ is even; the values for odd $d$ were given in 
Theorem \ref{T:odd3}.
\vskip 10pt

\begin{center}
\begin{tabular}{|c|c|c|}
\hline 
$disc(q)$ & $\epsilon_{HW}(q)$ & $\lambda_L$ \\
\hline 
1 & $(-1)^{n(n-1)/2}$ & $(2^{n-1}+1)(2^n-1)/2$ \\  
\hline 
1 & $-(-1)^{n(n-1)/2}$ & $(2^{n-1}-1)(2^n-1)/2$ \\
\hline   
$E_{disc(q)}$ unramified & $(-1)^{n(n-1)/2}$ & $(2^{n-1}-1)(2^n+1)/2$ \\  
\hline 
$E_{disc(q)}$ unramified & $-(-1)^{n(n-1)/2}$ & $(2^{n-1}+1)(2^n+1)/2$ \\
\hline   
$E_{disc(q)}$ ramified  & $\pm 1$ & $1/2$ \\
\hline
\end{tabular}
\end{center}  

\vskip 10pt

\section{\bf Mass of Unimodular Lattices}
Assembling the results of the previous section, and 
using Theorem \ref{T:mass},
one can give explicit formulas for the mass of 
unimodular lattices over $\mathbb{Z}$.
We shall only consider the odd unimodular lattices here, 
since the even case is easy. For the following proposition, see [Se].

\begin{Prop}
Let $r \geq s \geq 0$ with $r+s =d$.
Let $(V_{r,s},q_{r,s})$ be the quadratic space 
$\mathbb{Q}^d$ defined by the quadratic form
\[ q_{r,s}(x) = \sum_{i=1}^r x_i^2 - \sum_{j = r+1}^{r+s} x_j^2. \] 
Let $L_{r,s} = \mathbb{Z}^d$. 
\vskip 5pt

(i) Any odd unimodular lattice is in the genus of some $L_{r,s}$.
\vskip 5pt

(ii) If $s>0$, so that $V_{r,s}$ is indefinite over $\mathbb{R}$,
then any odd unimodular lattice is isomorphic to some $L_{r,s}$.
\end{Prop}

The proposition shows that there is no loss of generality
in working with the lattice 
$L = L_{r,s}$ above. Over a finite place $p$, 
the basic invariants of $V_{r,s}$ are
\[ disc(q_{r,s}) = (-1)^{[d/2]+ s} \quad \text{and} \quad 
\epsilon_{HW}(q_{r,s}) = (-1,-1)^{[s/2]}. \]
If $p \ne 2$, then $L_{r,s} \otimes \mathbb{Z}_p$ is a maximal lattice in the
relevant quadratic space and a quick check shows that the $\lambda$-factor is 
1. On the other hand, when $p=2$, we can read off the value of 
$\lambda_{L \otimes \mathbb{Z}_2}$ from the two tables of the previous section.
We state the results for odd and even dimensional spaces separately:

\begin{Thm} \label{T:mass2}
Assume that $d = 2n+1$ is odd. Then
\[ Mass(L_{r,s}) = \lambda_2 \cdot 
\prod_{k=1}^n \frac{(2k-1)! \cdot \zeta(2k)}{(2 \pi)^{2k}} \cdot \tau(G)  \]
where $\tau(G)=2$ and
\[ \lambda_2 = \begin{cases}
(2^n+1)/2, \text{   if $r-s \equiv \pm 1$ (mod 8);} \\
(2^n -1)/2,  \text{   if $r-s \equiv \pm 3$ (mod 8).} \end{cases} \]
\end{Thm}
\vskip 10pt

\begin{Thm} \label{T:mass3}
Assume that $d = 2n$ is even. Then
\[ Mass(L_{r,s}) = \lambda_2 \cdot 
\frac{(n-1)! \cdot L(n, \chi_{(-1)^{(r-s)/4}})}{(2 \pi)^n} \cdot  
\prod_{k=1}^{n-1} \frac{ (2k-1)! \cdot \zeta(2k)}{(2\pi)^{2k}} 
\cdot d_{\mathbb{Q}((-1)^{(r-s)/4})}^{n-1/2}  \cdot \tau(G)
 \]
where $\tau(G)=2$,
\[  
d_{\mathbb{Q}((-1)^{(r-s)/4})} = \begin{cases}
4 \text{  if $r-s \equiv \pm 2$ (mod 8);} \\
1 \text{  if $r-s \equiv 0$ or $4$ (mod 8),} \end{cases} \]
and
\[ \lambda_2 = \begin{cases}
1/2, \text{   if $r-s \equiv \pm 2$ (mod 8);} \\
(2^{n-1}+1)(2^n -1)/2,  \text{   if $r-s \equiv 0$ (mod 8);} \\ 
(2^{n-1}-1)(2^n -1)/2,  \text{   if $r-s \equiv 4$ (mod 8).}  
\end{cases} \]
\end{Thm}
\vskip 10pt

In the following tables, we give the values of masses for small $n$.
\bigskip

{\bf 1.} $r+s = 2n +1$ is odd:
$$
\def\arraystretch{1.5}
\begin{array}{|c|c|c|}
\hline
  n  & r-s \equiv \pm 1 (8) & r-s \equiv \pm 3 (8) \\
\hline
1 & \frac{1}{8} & \frac{1}{24} \\
2 & \frac{1}{1152} & \frac{1}{1920} \\
3 & \frac{1}{322560} & \frac{1}{414720} \\
4 & \frac{17}{1393459200} & \frac{1}{92897280}\\
5 & \frac{1}{11147673600} & \frac{31}{367873228800}\\
6 & \frac{691}{370816214630400} & \frac{691}{382588157952000}\\
7 & \frac{29713}{192824431607808000} & \frac{87757}{578473294823424000}\\
8 & \frac{642332179}{9440684171518279680000} & \frac{2499347}{37022290868699136000}\\
9 & \frac{109638854849}{528678313605023662080000} & \frac{8003636403977}{38744567839911019806720000}\\
10 & \frac{784910153445588299}{143199922736311129205637120000} & \frac{593468652605200909}{108484789951750855458816000000}\\
\hline
\end{array}
$$
\bigskip

{\bf 2.} $r+s = 2n$ is even:
$$
\def\arraystretch{1.5}
\begin{array}{|c|c|c|c|}
\hline
  n  & r-s \equiv 0(8) & r-s \equiv\pm 2(8) & r-s \equiv 4(8) \\
\hline
2 & \frac{1}{64} & \frac{L(2)}{12\;\pi^2} & \frac{1}{192} \\
3 & \frac{7\;\zeta(3)}{4608\;\pi^3} & \frac{1}{23040} & \frac{7\;\zeta(3)}{7680\;\pi^3} \\
4 & \frac{1}{5160960} & \frac{L(4)}{60480\;\pi^4} & \frac{1}{6635520} \\
5 & \frac{527\;\zeta(5)}{1857945600\;\pi^5} & \frac{1}{1114767360} & \frac{31\;\zeta(5)}{123863040\;\pi^5} \\
6 & \frac{1}{89181388800} & \frac{L(6)}{95800320\;\pi^6} & \frac{31}{2942985830400} \\
7 & \frac{87757\;\zeta(7)}{65922882600960\;\pi^7} & \frac{42151}{96412215803904000} & \frac{87757\;\zeta(7)}{68015672524800\;\pi^7} \\
8 & \frac{505121}{6170381811449856000} & \frac{691\;L(8)}{896690995200\;\pi^8} & \frac{1491869}{18511145434349568000} \\
9 & \frac{46890249067\;\zeta(9)}{17125957680758784000\;\pi^9} & \frac{692319119}{7552547337214623744000} & \frac{182452331\;\zeta(9)}{67160618355916800\;\pi^9} \\
10 & \frac{3398804500319}{4229426508840189296640000} & \frac{109638854849\;L(10)}{1459741204905984000\;\pi^{10}} & \frac{248112728523287}{309956542719288158453760000} \\
\hline
\end{array}
$$
\bigskip

\noindent{\bf Comments:}
The values of masses of the odd unimodular lattices first decrease when the dimension of the 
space grows and attain their minima $\sim 10^{-14}$ at $n = 8$ for each of the cases, then the 
mass starts to grow exponentially. All the values in the first table are 
rational numbers which is related to the fact
that the Euler-Poincar\'e characteristic for symmetric spaces of 
the corresponding orthogonal groups 
does not vanish. In the second table one can see many irrationalities 
some of which are of a particular
interest. Thus, for the type $(3,1)$ which corresponds to the hyperbolic 
$3$-space, we find that the value
of the mass is a rational multiple of Catalan's constant 
$$C = L(2) = 1 - \frac{1}{9} + \frac{1}{25}- \ldots $$
devided by $\pi^2$. It is conjectured but not known that both $C$ and $C/\pi^2$ are irrational.

\vskip 10pt
\noindent{\bf Acknowledgements:}
We would like to thank Professor Don Zagier for his comments and
suggestions. The first author would like to thank
the Max-Planck-Institut f\"ur Mathematik in Bonn
for hospitality and financial support.
The second author is partially supported by NSF grant DMS-0202989 and the
AMS Centennial Fellowship.

\end{document}